\documentclass[12pt]{article}

\usepackage{amsmath}
\usepackage{amssymb}
\usepackage{eucal}

\usepackage[english]{babel}

\begin{document}
 
\baselineskip 16pt
 
\title{Finite groups with $\Bbb P$-subnormal 2-maximal
subgroups}


\author{V.\,N. Kniahina and V.\,S. Monakhov}


\maketitle

\begin{abstract}
A subgroup $H$ of a group $G$ is called $\Bbb P$-{\sl subnormal} in $G$ 
if either $H=G$ or there is a chain of subgroups 
$H=H_0\subset H_1\subset \ldots \subset H_n=G$ 
such that $|H_i:H_{i-1}|$~is prime for $1\le i\le n$.
In this paper we study the groups all of whose
2-maximal subgroups are $\Bbb P$-subnormal.
\end{abstract}


MSC2010 20D20, 20E34


{\small }

\section{Introduction}

We consider finite groups only. A subgroup $K$ of a group $G$ is 
called 2-maximal in $G$ if $K$ is a maximal subgroup of some
maximal subgroup $M$ of $G$.

Let $H$ be a subgroup of a group $G$ and $n$ is a positive integer.
If there is a chain of subgroups
$$
H=H_0\subset H_1\subset \ldots \subset H_{n-1}\subset H_n=G,
$$ 
such that $H_i$ is a maximal subgroup of $H_{i+1}$, 
$i=0,1,\ldots ,n-1$, then $H$ is called $n$-maximal in $G$.

For example, in the symmetric group $S_4$ the subgroup $I$ of order 2 
from $S_3$ is 2-maximal in the chain of subgroups $I\subset S_3\subset S_4$ 
and 3-maximal in the chain of subgroups  
$I\subset Z_4\subset D_8\subset S_4$. Here, $Z_4$ is the cyclic group 
of order 4 and $D_8$ is the dihedral group of order 8.
For any $n\ge 3$, there exists a group in which some 
2\nobreakdash-\hspace{0pt}maximal subgroup is $n$-maximal, see 
Example 1 below.

A.\,F. Vasilyev, T.\,I. Vasilyeva and V.\,N. Tyutyanov in \cite{VVTGGU}
introduced the following definition. Let $\Bbb P$ be the set of all prime 
numbers. A subgroup $H$ of a group $G$ is called $\Bbb P$-{\sl subnormal} in $G$ 
if either $H=G$ or there is a chain  
$$
H=H_0\subset H_1\subset \ldots \subset H_n=G
$$
of subgroups such that $|H_i:H_{i-1}|$~is prime for  
$1\le i\le n$. 
In \cite{VVTGGU}, \cite{VVTSMJ} studied groups with $\Bbb P$-subnormal 
Sylow subgroups. 

In \cite{VVTGGU} proposed the following problem:

\medskip

{\sl Describe the groups in which all 2-maximal subgroups are 
$\Bbb P$\nobreakdash-\hspace{0pt}subnormal.}

\medskip

This problem is solved in the article. The following 
theorem is proved. 

\medskip

{\bf Theorem.} {\sl Every 2-maximal subgroup of a group $G$
is $\Bbb P$\nobreakdash-\hspace{0pt}subnormal in $G$
if and only if~${\Phi (G^{\frak U})=1}$ and every proper 
subgroup of $G$ is supersolvable.}

\medskip

Here, $G^{\frak U}$ is the smallest normal subgroup of $G$ such that the
corresponding quotent group is supersolvable, $\Phi (G^{\frak U})$~is the 
Frattini subgroup of $G^{\frak U}$.

\section{Preliminary results}

We use the standart notation of \cite{Hup}. The set of prime divisors
of $|G|$ is denoted $\pi (G)$. We write $[A]B$ for a semidirect product 
with a normal subgroup $A$. If $H$ is a subgroup of a 
group $G$, then $\bigcap _{x\in G}x^{-1}Hx$ is called the core of $H$~in~$G$,
denoted $H_G$. If a group $G$ contains a maximal subgroup $M$ with trivial
core, then $G$ is said to be primitive and $M$ is its primitivator.
We will use the following notation: $S_n$ and $A_n$ are the symmetric and 
alternating groups of degree $n$, $E_{p^t}$~is the elementary abelian group
of order $p^t$, $Z_m$~is the cyclic group of order $m$. Let 
$|G|=p_1^{a_1}p_2^{a_2} \ldots p_k^{a_k}$, where
$p_1>p_2> \ldots >p_k$. We say that $G$ has an ordered Sylow tower
of supersolvable type if there exist normal subgroups $G_i$ with
$$
1=G_0\leq G_1\leq G_2\leq \ldots \leq G_{k-1}\leq G_k=G,
$$ 
and, where each factor $G_{i}/G_{i-1}$ is 
isomorphic to a Sylow $p_i$-subgroup of $G$ for all $i=1,2,\ldots,k$.

\medskip

{\bf Lemma 1.} \cite[Theorem IX.8.3]{HB2}
{\sl  Let $a$, $n$ be integers greater than 1. Then except in the cases
$n=2$, $a=2^b-1$ and $n=6$, $a=2$, there is a prime $q$ with the 
following properties:
                                                    
$1)$ $q$ divides $a^n-1$;

$2)$ $q$ does not divide $a^i-1$ whenever $0<i<n$;

$3)$ $q$ does not divide $n$.}

\medskip  

{\bf Example 1.}                                 
For every $n\ge 3$ there exists a group in which some 2-maximal subgroup is
$n$-maximal. Let $n=3$. In the symmetric group $S_4$  
the subgroup $I$ of order 2 from $S_3$ is  
2-maximal in the chain of subgroups $I\subset S_3\subset S_4$ 
and 3-maximal in the chain of subgroups  
$I\subset E_4\subset D_8\subset S_4$. Now let $n>3$ and $a=5$. By  
Lemma 1, there exists  
a prime $q$ such that $q$ divides $5^{n-1}-1$ and $q$ does  not divide 
$5^i-1$ for all $i\in \{1,2,\ldots ,n-2\}$. Hence $GL(n-1,5)$
contains a subgroup $Z$ of order $q$ which acts irreducibly on the 
elementary abelian group $E_{5^{n-1}}$ of order $5^{n-1}$.
In the group $X=[E_{5^{n-1}}]Z$ the identity  subgroup 1 is 2-maximal  
in the chain of subgroups $1\subset Z\subset X$ and $n$-maximal
in the chain of subgroups  
$1\subset E_5\subset E_{5^2}\subset \ldots \subset E_{5^{n-1}}\subset X$.

\medskip 

Recall that a Schmidt group is a finite non-nilpotent group in which every 
proper subgroup is nilpotent.
 
\medskip 

{\bf Example 2.}
Let $S=[P]Q$ be a Schmidt group of order $2^{11}11$, $A=\Phi (P)$,
$|A|=2$. Then $A\times Q$ is maximal in $S$, $A$ is 2-maximal in $S$,
and $A$ is 10-maximal in $S$ because 
$A=A_0\subset A_1\subset \ldots \subset A_{9}=P\subset S$, 
$|A_i:A_{i-1}|=2$, $1\le i\le 9$, $|S:P|=11$.

\medskip

{\bf Lemma 2.} \cite[Lemma 2.1]{VVTGGU} {\sl 
Let $N$~be a normal subgroup of a group $G$, $H$~an
arbitrary subgroup of $G$. Then the following hold:

$1)$ if $H$ is $\Bbb P$-subnormal in $G$, then 
$(H\cap N)$ is $\Bbb P$-subnormal in $N$, and $HN/N$ is
$\Bbb P$-subnormal in $G/N$;

$2)$ if $N\subseteq H$ and $H/N$ is $\Bbb P$-subnormal in $G/N$, 
then $H$ is $\Bbb P$-subnormal in $G$;

$3)$ if $H$ is $\Bbb P$-subnormal in $K$, $K$ is $\Bbb P$-subnormal in 
$G$, then $H$ $\Bbb P$-subnormal in $G$;

$4)$ if $H$ is $\Bbb P$-subnormal in $G$, then $H^g$ is $\Bbb P$-subnormal in 
$G$ for each element $g\in G$}.

\medskip 

{\bf Example 3.}
In the alternating group $G=A_5$ the subgroup $H=A_4$ is $\Bbb P$-subnormal.
If $x\in G\setminus H$, then $H^x$ is $\Bbb P$-subnormal in $G$.
The subgroup $D=H\cap H^x$ is a Sylow 3-subgroup of the group $G$ and $D$
is not $\Bbb P$\nobreakdash-\hspace{0pt}subnormal in $H$. 
Therefore an intersection of two $\Bbb P$-subnormal subgroups is not 
$\Bbb P$-subnormal.
Moreover, if a subgroup $H$ is $\Bbb P$-subnormal in a group $G$ and 
$K$~is an arbitrary subgroup of $G$, in general, their intersection 
$H\cap K$ is not $\Bbb P$-subnormal in~$K$. 

\medskip

{\bf Lemma 3.}
{\sl Let $H$ be a subgroup of a solvable group $G$, and 
assume that $|G:H|$~is a prime number. Then $G/H_G$ is supersolvale.}

\medskip

{\sc Proof.}
By hypothesis, $|G:H|=p$, where $p$ is a prime number. If $H=H_G$, then $G/H$~is
cyclic of  prime order $p$, and thus $G/H_G$ is supersolvable, 
as required.  Assume now that $H\ne H_G$, i.\,e.
$H$ is not normal in $G$. Then $G/H_G$ contains a maximal subgroup
$H/H_G$ with trivial core. Therefore $G/H_G$ is primitive and 
its Fitting subgroup $F/H_G$ has prime order $p$. 
Since $F/H_G=C_{G/H_G}(F/H_G)$, it follows that 
$(G/H_G)/(F/H_G)\simeq H/H_G$ is isomorphic to a cyclic group 
of order dividing $p-1$. Thus $G/H_G$ is supersolvable.

\medskip                             

{\bf Lemma 4.} {\sl Let $p$ be the largest prime divisor of  
$|G|$, and suppose that $P$ is a Sylow $p$\nobreakdash-\hspace{0pt}subgroup
of $G$. Assume that  
$P$ is not normal in $G$, and that $H,K\subseteq G$ are subgroups with 
$N_G(P)\subseteq K\subseteq H$. Then $|H:K|$ is not prime.}

\medskip

{\sc Proof.} It is clear that $N_G(P)=N_K(P)=N_H(P)$, and $P$ is a Sylow 
$p$-subgroup of $K$ and of $H$. By the lemma on indexes, we have
 
$$
|H:N_H(P)|=|H:K|\ |K:N_K(P)|,
$$
and, by the Sylow theorem, 
$$
|H:N_H(P)|=1+hp, \ |K:N_K(P)|=1+kp, \ h, \ k\in \Bbb N\cup \{0\}.
$$
Let $|H:K|=t$. Now, 
$$
1+hp=t(1+kp), \ \ t=1+(h-tk)p.
$$
We see that $p$ divides $t-1$, and thus $t>p$. 
If $t$ is prime, this contradicts the maximality of $p$.

{\bf Lemma 5.}  1. {\sl A group is supersolvable if and only if
the index of every of its maximal subgroup is prime.

$2$.  Every subgroup of a supersolvable group is $\Bbb P$-subnormal.

$3$.  A group is supersolvable if and only if the normalizers of all
of its Sylow subgroups are $\Bbb P$-subnormal.}

\medskip

{\sc Proof.} 1. This is Huppert's classic 
theorem, see \cite[Theorem VI.9.5]{Hup}.

2. The statement follows from (1) of the 
lemma.

3.  If a group is supersolvable, then all of its subgroups are 
$\Bbb P$\nobreakdash-\hspace{0pt}subnormal, see (2).

Conversely, suppose that the normalizer of every Sylow subgroup
of a group $G$ is $\Bbb P$-subnormal. By 
Lemma 4, for the largest 
$p\in \pi (G)$ a Sylow $p$-subgroup $P$ of $G$ is normal. 
It is easy to check that the conditions of the lemma are
inherited by all quotient groups and so $G/P$ is supersolvable.
In particular, $G$ has an ordered Sylow tower of supersolvable type. 
Since the class of all supersolvable groups is a saturated formation, 
we can assume, by the inductive hypothesis, that $G$ is primitive, in particular, 
$G=[P]M$, where $M$~is a maximal subgroup with trivial core.
Since $M$ is supersolvable, it follows that $M=N_G(Q)$ for the largest 
$q\in \pi (M)$.
It is obvious that $p\ne q$ and $M=N_G(Q)$ is $\Bbb P$-subnormal in $G$, by
the condition of the lemma. Therefore $|P|=p$ and, by 
Lemma 3, $G$ is supersolvable.

\medskip

{\bf Lemma 6.} {\cite[Theorem 22]{Hup54}, \cite{Doerk} {\sl
Let $G$~be a minimal non-supersolvable group. We have:

$1)$ $G$ is solvable and $|\pi (G)|\le 3$.

$2)$ If $G$ is not a Schmidt group, then $G$ has an ordered 
Sylow tower of supersolvable type.

$3)$ $G$ has a unique normal Sylow subgroup $P$
and $P=G^{\frak U}$.

$4)$ $|P/\Phi(P)|>p$ and $P/\Phi (P)$ is a minimal normal subgroup    
of~$G/\Phi (G)$.

$5)$ The Frattini subgroup $\Phi (P)$ of $P$ is supersolvable 
embedded in $G$, i.e., there exists a series
$$
1\subset N_0\subset N_1 \ldots \subset  N_n=\Phi (P)
$$
such that $N_i$ is a normal subgroup of $G$ and $|N_i/N_{i-1}|\in \Bbb P$ for
$1\le i\le n$.

$6)$ Let $Q$ be a complement to $P$ in $G$. Then $Q/Q\cap \Phi (G)$~is
a minimal non-abelian group or a cyclic group of prime power order.

$7)$ All maximal subgroups of non-prime index are conjugate in $G$,
and moreover, they are conjugate to $\Phi (P)Q$.}

\bigskip

\section{Main results}

\medskip

{\bf Theorem.} {\sl Every 2-maximal subgroup of a group $G$
is $\Bbb P$\nobreakdash-\hspace{0pt}subnormal in $G$
if and only if ${\Phi (G^{\frak U})=1}$ and every proper 
subgroup of $G$ is supersolvable.}

\medskip

{\sc Proof.} Suppose that all 2-maximal subgroups of a group
$G$ are $\Bbb P$-subnormal. We proceed by induction on $|G|$. Show
first that $G$ has an ordered Sylow tower of supersolvable type.
By Lemma 2, the conditions of the theorem are inherited by all quotient    
groups of $G$. 

\medskip

(1) $G$ has an ordered Sylow tower of supersolvable type.

Let $P$~be a Sylow $p$\nobreakdash-\hspace{0pt}subgroup of $G$, 
where
$p$ is the largest prime divisor of $|G|$.  
Suppose that $P$ is not normal in $G$. It follows that $N_G(P)$ is 
a proper subgroup of $G$. If $N_G(P)$ is not maximal in $G$, then 
there exists a 2-maximal subgroup $A$ containing $N_G(P)$. 
By the condition of the theorem, $A$ is $\Bbb P$\nobreakdash-\hspace{0pt}
subnormal in $G$, and so $A$ is contained in a subgroup of prime index. 
This contradicts Lemma~4. Therefore $N_G(P)$ is maximal in $G$ 
and $|G:N_G(P)|\not \in \Bbb P$ by Lemma~4. 
If $N_G(P)=P$, then $G$ is solvable by Theorem IV.7.4 \cite{Hup}. 
It follows that $N_G(P)\ne P$ and $N_G(P)$ has a maximal subgroup 
$B$ which contains $P$. We see that $B$ is 2-maximal in $G$ and, 
by the condition
of the theorem, $B$ is $\Bbb P$\nobreakdash-\hspace{0pt}subnormal. 
Hence there exists a chain of subgroups 
$$
P\subseteq B=B_0\subset  B_1\subset \ldots \subset B_{t-1}=V\subset  B_t=G,
\ |B_i:B_{i-1}|\in \Bbb P, \ 1\le i\le t.
$$

The subgroup $V$ is maximal in $G$ and $V$ different from $N_G(P)$, because 
$|G:N_G(P)|$~is not a prime number, whereas $|G:V|$~is prime. Besides,
$t\ge 3$. Thus $V\cap N_G(P)=B$ and $N_V(P)=V\cap N_G(P)=B=N_{B_1}(P)$. 
We have $|B_1:N_{B_1}(P)|\in \Bbb P$, this contradicts 
Lemma~4. Therefore the assumption is false and $P$ is normal in $G$.
By induction on $|G|$, every proper subgroup of $G/P$ is supersolvable, 
and by Lemma 6, $G/P$ has an ordered Sylow tower of supersolvable type.
Thus $G$ has an ordered Sylow tower of supersolvable type, 
in particular, $G$ is solvable.

\medskip

(2) Every proper subgroup of $G$ is supersolvable.

Suppose that $G$ contains a non-supersolvable maximal subgroup $H$. 
Then, by Lemma 5, $H$ contains a maximal subgroup $K$ of 
non-prime index.
Since $K$ is 2-maximal in $G$, there exists a chain of subgroups 
$$
K=K_0\subset K_1\subset \ldots \subset  K_{n-1}=T\subset K_n=G
$$
such that $|K_i:K_{i-1}|\in \Bbb P$ for all $i=1,2,\ldots , n$. 
It is clear that $H\ne T$ and $H\cap T=K$.

Assume that $G=HT$. In this case, 
$$
|G:T|=|H:H\cap T|=|H:K|\in \Bbb P,
$$
this is a contradiction. Hence $G\ne HT$. Since $H$ and $T$~are 
distinct maximal subgroups of $G$, and $G$ is solvable, by Theorem
II.3.9 \cite{Hup}, we have $T=H^g$ for some $g\in G$.
Since $H\ne T$, we see that $H$~is a non-normal  maximal subgroup
of prime index in $G$. By Lemma 3, the quotient group $G/H_G$ 
is supersolvable. Since
$$
H_G\subseteq H\cap H^g=H\cap T=K,
$$
we have $K/H_G$ is maximal in $H/H_G$. 
By Lemma 5,
$$
|H:K|=|H/H_G:K/H_G|\in \Bbb P,
$$
this is a contradiction. Therefore the assumption is false and
every proper subgroup of $G$ is supersolvable.

\medskip

(3) ${\Phi (G^{\frak U})=1}$ 

If $G$ is supersolvavle, then $G^{\frak U}=1$, it follows that 
$\Phi (G^{\frak U})=1$. Assume now that $G$ is non-supersolvable. 
Then $G$ has the properties listed in Lemma 6. We keep the notation 
of that lemma. Now $G^{\frak U}=P$ and $[\Phi (P)]Q$~is maximal in $G$. 

Suppose that $\Phi (P)\ne 1$. Assume that $A=N_{m-1}$~is a maximal
subgroup of $\Phi (P)$, and that $A$ is normal in $G$. 
Then $[A]Q$~is a 2-maximal subgroup of $G$. By the condition of the 
theorem, 
$[A]Q$ is $\Bbb P$-subnormal in $G$. Hence, there exists a chain
of subgroups $[A]Q \subseteq B\subseteq G$ such that $|G:B|\in \Bbb P$.
Since $G=[P]Q$ and $Q\subseteq B$, by the Dedekind identity, we have 
$B=(B\cap P)Q$, and $B\cap P$~is maximal in $P$. Therefore
$\Phi (P)\subseteq B\cap P$ and $\Phi (P)Q$ is conained in $B$, where
$\Phi (P)Q$ is maximal in $G$.  
Thus $B=\Phi (P)Q$ and $p=|G:B|=|P:\Phi (P)|$, this contradicts 
Lemma 6. Therefore our assumption is false and $\Phi (P)=1$. 
The necessity is proved.

Prove the sufficiency. Assume that every proper subgroup of $G$ is
supersolvable and $\Phi (G^{\frak U})=1$. If a group is 
supersolvable, then every its maximal subgroup has a prime index,
it follows that every 2-maximal subgroup of a supersolvable group is 
$\Bbb P$\nobreakdash-\hspace{0pt}subnormal. 
Let $G$ be non-supersolvable. Then $G$ is minimal non-supersolvable
and the structure of $G$ is described in Lemma 6.
We keep for $G$ the notation of that lemma, in particular, we have:  
$P=G^{\frak U}$, $\Phi (P)=1$ and $Q$~is a maximal subgroup of $G$. 
Let $H$~be an arbitrary 2-maximal subgroup of the group $G$. If $H\subseteq M$,
where $M$ is a maximal subgroup of $G$ and $|G:M|\in \Bbb P$, then $H$ is
$\Bbb P$-subnormal in $G$, because $M$ is supersolvable.
If $H\subseteq K$, where $K$ is a maximal subgroup of the group $G$ 
and $|G:K|\not \in\Bbb P$, then, by Lemma 6, the subgroup $H$  
contained in $Q^g$ for some $g\in G$. Therefore $PH$~is a proper subgroup
of $G$, thus $PH$ is supersolvable, and $H$~is $\Bbb P$-subnormal in $PH$. 
Let $T$~be a maximal subgroup of $G$ containing $PH$. Since $T$
is supersolvable and $|G:T|\in \Bbb P$, we see that 
$PH$~is $\Bbb P$-subnormal in $G$. Using Lemma 2, we deduce that
$H$~is $\Bbb P$-subnormal in $G$. The theorem is proved.

\medskip

{\bf Corollary.} {\sl Suppose that every 2-maximal subgroup of a group 
$G$ is $\Bbb P$-subnormal. If $|\pi (G)|\ge 4$, then $G$ is supersolvable.}

\medskip                                     

{\sc Proof.} Let every 2-maximal subgroup of a group $G$ be $\Bbb P$-subnormal. 
Suppose that $G$ is not supersolvable. By the previous 
theorem, $G$ is a minimal non-supersolvable group. By Lemma 6, the order of $G$ 
has at most three prime divisors, i.e. $|\pi (G)|\le 3$, which is a contradiction. 
Therefore, our assumption is false and $G$ is supersolvable.

\medskip
 
The following examples show that for $|\pi (G)|=2$ and for $|\pi (G)|=3$
there exist non-supersolvable groups in which every 2-maximal subgroup 
is $\Bbb P$-subnormal.       
\medskip 

{\bf Example 4.}
There are three non-isomorphic minimal non-supersolvable groups of order 400: 
$$
[E_{5^2}](<a><b>), \ |a|=|b|=4.
$$
Numbers of these groups in the library of SmallGroups~\cite{GAP}
are [400,129], 
[400,130], [400,134]. The Sylow 2-subgroups of these groups are 
non-abelian and have the form: $[Z_4\times Z_2]Z_2$ and $[Z_4]Z_4$. 
Suppose that $G$ is one of these groups.
Then $G^{\frak U}=[E_{5^2}]$ and $\Phi (G^{\frak U})=1$. All subgroups
of the group $G$ are $\Bbb P$-subnormal, except the maximal subgroup 
$<a><b>$. 

\medskip 

{\bf Example 5.}
The general linear group $GL(2,7)$ contains the symmetric group $S_3$
which acts irreducibly on the elementary abelian group $E_{7^2}$ of order 49. 
The semidirect product $[E_{7^2}]S_3$ is a minimal non-supersolvable 
group, it has subgroups of orders 14 and 21.
Therefore, in the group $[E_{7^2}]S_3$, every 2-maximal subgroup is 
$\Bbb P$\nobreakdash-\hspace{0pt}subnormal.

\bigskip


\noindent V.\,N. KNIAHINA

\noindent  Gomel Engineering Institute, Gomel
246035, BELARUS

\noindent E-mail address: knyagina@inbox.ru

\bigskip

\noindent V.\,S. MONAKHOV

\noindent Department of mathematics, Gomel F. Scorina State
University, Gomel 246019,  BELARUS 

\noindent E-mail address: Victor.Monakhov@gmail.com

\end{document}